\documentclass{article}
\newtheorem{thm}{Theorem}[section]
\newtheorem{cor}[thm]{Corollary}
\newtheorem{prop}[thm]{Proposition}

\newtheorem{lem}[thm]{Lemma}
\newtheorem{Def}[thm]{Definition}

\newtheorem{ex}[thm]{Example}

\newcommand{\be}{\begin{equation}}
\newcommand{\ee}{\end{equation}}
\newcommand{\ben}{\begin{enumerate}}
\newcommand{\een}{\end{enumerate}}
\newcommand{\beq}{\begin{eqnarray}}
\newcommand{\eeq}{\end{eqnarray}}
\newcommand{\beqn}{\begin{eqnarray*}}
\newcommand{\eeqn}{\end{eqnarray*}}
\newcommand{\e}{\varepsilon}
\newcommand{\pa}{\partial}

\newcommand{\pxi}{ {\pa \over \pa x^i}}

\newcommand{\pyi}{{\pa \over \pa y^i}}

\newcommand{\qed}{\hspace*{\fill}Q.E.D.} 

\title{A Class of Isotropic Mean Berwald Metrics}

\author{B. Najafi and A. Tayebi}

\begin{document}

\maketitle
\begin{abstract}
In this paper, we find a condition on $(\alpha, \beta)$-metrics  under which the notions of isotropic S-curvature, weakly isotropic S-curvature and isotropic mean Berwald curvature are equivalent.\footnote{2010 {\it Mathematics Subject Classification}:  Primary
53B40, 53C60}\\\\
{\bf {Keywords}}: $(\alpha, \beta)$-metric, isotropic S-curvature,  isotropic E-curvature.
\end{abstract}

\section{Introduction}
The S-curvature is introduced  by Shen for a comparison
theorem on Finsler manifolds \cite{Sh3}. Recent study shows that the S-curvature
plays a very important role in Finsler geometry \cite{SX}\cite{TR}. A  Finsler metric $F$ is said to have  isotropic S-curvature if ${\bf S}= (n+1) c F$,  where $c=c(x)$ is a  scalar function on an n-dimensional manifold $M$.

Taking twice vertical covariant derivatives  of the S-curvature gives rise the $E$-curvature. A Finsler metric $F$ with vanishing $E$-curvature called weakly Berwald metric. In \cite{BY}, B\'{a}cs\'{o}-Yoshikawa study some weakly Berwald metrics. Also, $F$  is called to have isotropic
E-curvature if ${\bf E}= \frac{n+1}{2} c F^{-1}{\bf h}$,  for some scalar function
$c$ on $M$, where ${\bf h}$ is the angular metric. It is easy to see that every Finsler metric of isotropic
S-curvature is an isotropic E-curvature.  Now, is the equation
${\bf S}= (n+1) c F$ equivalent to the equation ${\bf E}= \frac{n+1}{2} c F^{-1}{\bf h}$?

Recently, Cheng-Shen prove that a Randers metric $F=\alpha+\beta$ is of isotropic $S$-curvature if and only if it is of isotropic $E$-curvature \cite{ChSh2}. Then,  Chun-Huan-Cheng extend this equivalency to the Finsler metric  $F=\alpha^{-m}(\alpha+\beta)^{m+1}$ for every real constant $m$,  including Randers metric \cite{CC}. In \cite{LL}, Lee-Lee prove that  this notions are  equivalent for the Finsler metrics in the form  $F=\alpha+\alpha^{-1}\beta^2$.

 All of above metrics are special Finsler metrics so- called $(\alpha,\beta)$-metrics.  An $(\alpha, \beta)$-metric is a scalar function on $TM$ defined by $F:=\alpha\phi(s)$ , $s=\beta/\alpha$ where  $\phi=\phi(s)$ is a $C^\infty$ on $(-b_0, b_0)$ with certain regularity, $\alpha$ is a Riemannian metric and  $\beta$  is a 1-form on a manifold $M$ \cite{TP}.  A natural question arises:\\

{\it Is being of isotropic $S$-curvature is equivalent to being of isotropic $E$-curvature for $(\alpha,\beta)$-metrics?}\\

In  \cite{DW},  Deng-Wang find the formula of the $S$-curvature of homogeneous $(\alpha,\beta)$-metrics. Then  Cheng-Shen classify $(\alpha,\beta)$-metrics of isotropic $S$-curvature \cite{ChSh3}.

Let $F=\alpha\phi(s)$ be an $(\alpha,\beta)$-metric on a manifold $M$ of dimension n, where $s=\frac{\beta}{\alpha}$, $\alpha=\sqrt{a_{ij}y^iy^j}$ is a Riemannian metric and  $\beta=b_i(x)y^i$ is a 1-form on $M$.  For an $(\alpha,\beta)$-metric, put
\begin{eqnarray*}
Q\!\!\!\!&=&\!\!\!\! { \phi'\over  \phi - s \phi' },\\
\Delta\!\!\!\!&=&\!\!\!\!  1+sQ + (b^2-s^2) Q',\\
\Phi\!\!\!\!&=&\!\!\!\! - (Q -s Q') \{ n\Delta + 1 + sQ \} -(b^2-s^2) (1+sQ) Q'',\\
\end{eqnarray*}
Using the same method as in \cite{ChSh3}, we give an affirmative answer to the above question for almost all $(\alpha,\beta)$-metrics. More precisely, we prove the following.
\begin{thm}\label{mainthm}
Let $F=\alpha\phi(s)$,  $s={\beta}/{\alpha}$, be an  $(\alpha,\beta)$-metric on a manifold $M$.  Let us define
\be
\Xi:=\frac{(b^2 Q+s) \Phi}{\Delta^2}.
\ee
Suppose that $\Xi$ is not constant. Then $F$ is of isotropic $S$-curvature if and only if it is of isotropic $E$-curvature.
\end{thm}

It is remarkable that if $\Xi=0$, then $F$ reduces to a Riemannian metric. But, in general, it is still an open problem if Theorem \ref{mainthm}
is true when $\Xi$ is a constant.
\begin{ex}
The above mentioned $(\alpha,\beta)$-metric correspond to  $\phi=1+s$, $\phi=(1+s)^{m+1}$ and $\phi=1+s^2$, respectively.
Using a Maple program shows that for all these metrics $\Xi$ is not constant.
\end{ex}

\section{Preliminaries}
Let $F=F(x,y)$ be a Finsler metric on an $n$-dimensional manifold $M$.
  There is a notion of distortion $\tau =\tau(x,y) $ on $TM$   associated with a volume form $dV=\sigma(x) dx$, which is defined by
\[
\tau  (x,y) = \ln \frac{ \sqrt{ \det (g_{ij}(x,y)) }}{ \sigma(x) }.
\]
Then the S-curvature is defined by
\[ {\bf S}(x,y) = \frac{d}{dt} \Big [ \tau \Big (c(t), \dot{c}(t) \Big ) \Big ] |_{t=0},\]
where $ c(t)$ is the geodesic with $c(0)=x$ and $\dot{c}(0)=y$ \cite{ChernShen}.
From the definition, we see that the S-curvature ${\bf S}(y)$ measures the rate of change in the distortion on $(T_{x}M, F_{x})$ in the direction $y\in T_{x}M$.

Let $\textbf{G}=y^i\pxi-2 G^i\pyi$ denote the spray of $F$ and $dV_{BH}=\sigma(x) dx$ be the Busemann-Hausdorff volume form on $M$, where the spray coefficients $G^i$ are defined by
\[
G^i(y):=\frac{1}{4}g^{il}(y)\Big\{\frac{\partial^2[F^2]}{\partial x^k
\partial y^l}(y)y^k-\frac{\partial[F^2]}{\partial x^l}(y)\Big\},\ \
y\in T_xM.
\]
Then the S-curvature is given by
\[
{\bf S}= \frac{\pa G^m}{\pa y^m} - y^m \frac{\pa }{\pa x^m} ( \ln \sigma ).
\]
The E-curvature ${\bf E}= E_{ij} dx^i\otimes dx^j$  is given by
\[ E_{ij} = \frac{1}{2} \frac{\pa^2 S }{\pa y^i\pa y^j}.\]
\begin{Def}\emph{
Let $(M, F)$ be a n-dimensional  Finsler manifold.  Then
\begin{description}
  \item[(a)]  $F$ is of isotropic S-curvature if ${\bf S}=(n+1)cF$,
  \item[(b)]  $F$ is of weak isotropic S-curvature if ${\bf S}= (n+1) c F+\eta$,
  \item[(c)]  $F$ is of isotropic E-curvature if ${\bf E}= \frac{n+1}{2} c F^{-1}{\bf h}$,
   \end{description}
where $c= c(x)$ is a scalar function on $M$,  $\eta=\eta_i(x)y^i$ is a $1$-form on $M$ and ${\bf h}$ is the angular metric \cite{ShDiff}.}
\end{Def}
Consider the $(\alpha,\beta)$-metric $F=\alpha\phi\Big (\frac{\beta}{\alpha}\Big)$ where  $\alpha=\sqrt{a_{ij}(x)y^iy^j}$ is a Riemannian metric and  $\beta=b_i(x)y^i$ is a $1$-form on a manifold $M$. In \cite{ChSh3}, the following is proved:

\begin{lem} Let $(M, F)$ be a n-dimensional Finsler manifold. Suppose that $F=\alpha \phi (\beta/\alpha)$ be an $(\alpha,\beta)$-metric. Then
\be \label{Gm_m}
\frac{\pa G^m}{\pa y^m} = y^m \frac{\pa }{\pa x^m} (\ln \sigma_{\alpha} )
+ 2 \Psi (r_0+s_0)-\alpha^{-1} \frac{\Phi}{2\Delta^2} (r_{00} -2\alpha Q s_0 ),
\ee
where
\begin{eqnarray*}
Q & = & \frac{\phi'}{\phi-s\phi'},\,\,\, \Delta  =  1+sQ + (b^2-s^2) Q',\,\,\,\Psi  = \frac{Q'}{2\Delta}\\
\Phi & = & - (Q -s Q') \{ n\Delta + 1 + sQ \} -(b^2-s^2) (1+sQ) Q''.
\end{eqnarray*}
\end{lem}
For an $(\alpha,\beta)$-metric, put
\begin{eqnarray*}
r_{ij}:= {1\over 2}  ( b_{i|j}+b_{j|i} )\!\!\!\!&,&\!\!\!\!\  s_{ij} := {1\over 2} ( b_{i|j} - b_{j|i}),\\
r_j := b^i r_{ij}, \   s_j:=b^i s_{ij},\   r_{i0}: = r_{ij}y^j \!\!\!\!&,&\!\!\!\!\ \  s_{i0}:= s_{ij}y^j,\  r_0:= r_j y^j,\  s_0 := s_j y^j.
\end{eqnarray*}
Let $\bar{G}^i$ denote the spray coefficients of $\alpha$ . We have the following formula for the spray coefficients $G^i$ of $F$ \cite{ChernShen}:
\[
G^i = \bar{G}^i + \alpha Q s^i_{\ 0} + \Theta \Big \{ -2 Q \alpha s_0 + r_{00} \Big \} \frac{y^i}{\alpha} + \Psi \Big \{ -2 Q \alpha s_0 + r_{00} \Big \}  b^i,
\]
where  $ s^{i}_{\ j}:=a^{ih}s_{hj},  s^{i}_{\ 0}:=s^{i}_{\ j}y^{j}$ and $r_{00}:=r_{ij}y^{i}y^{j}$. In \cite{ChSh3}, Cheng-Shen find   the S-curvature as follows
\be
{\bf S}
= \Big \{ 2\Psi  - \frac{f'(b)}{bf(b)}\Big \} (r_0+s_0)   -  \alpha^{-1}\frac{\Phi}{2\Delta^2}( r_{00}-2\alpha Q s_0) .  \label{S0}
\ee
Recently, Cheng-Shen  characterize $(\alpha,\beta)$-metrics with isotropic
S-curvature and proved the following.

\begin{lem} \label{Cheng-Shen}{\rm(\cite{ChSh3})} Let  $ F =\alpha \phi (\beta/\alpha)$ be an $(\alpha,\beta)$-metric on an $n$-manifold. The, $F$ is of isotropic S-curvature $ {\bf S} = (n+1) cF$, if and only if one of the following holds
\item[(i)] $\beta $ satisfies
\be
r_{ij} = \e \Big \{ b^2 a_{ij} -b_i b_j \Big \}, \ \ \ \ \ s_j=0, \label{rrrr1}
\ee
where $\e= \e(x)$ is a scalar function, and $\phi=\phi(s)$ satisfies
\be
 \Phi = -2 (n+1) k \frac{ \phi \Delta^2}{b^2-s^2}, \label{ode_phi}
\ee
where $k$ is a constant. In this case, $c=k\epsilon$.
\item[(ii)] $\beta$ satisfies
\be
r_{ij} =0,\ \ \ \ \ s_j =0.\label{rrrr2}
\ee
In this case, $ c =0$.
\end{lem}
It is remarkable that,   Cheng-Wang-Wang prove that the condition $\Phi=0$ characterizes the Riemannian metrics among $(\alpha,\beta)$-metrics \cite{ChWW}. Hence, in the continue, we suppose that $\Phi\neq 0$.
\section{Proof of Theorem \ref{mainthm}}
First, we find the formula  of E-curvature of $(\alpha,\beta)$-metrics. After a long and tedious computations, we obtain the following.
\begin{prop} \label{E}
Let $F=\alpha\phi({\frac{\beta}{\alpha}})$ be an $(\alpha,\beta)$-metric. Put $\Omega:=\frac{\Phi}{2\Delta^2}$. Then the E-curvature of  $F$ is given by the following
\begin{eqnarray}
\nonumber E_{ij}\!\!\!\!&=&\!\!\!\!\ C_1b_ib_j+C_2(b_iy_j+b_jy_i)+C_3y_iy_j+C_4a_{ij}+C_5(r_{i0}b_j+r_{j0}b_i)\\
\!\!\!\!&&\!\!\!\!\ \nonumber +C_6(r_{i0}y_j+r_{j0}y_i)+C_7r_{ij}+C_8(s_ib_j+s_jb_i)+C_9(s_iy_j+s_jy_i) \\
\!\!\!\!&&\!\!\!\!\ +C_{10}(r_ib_j+r_jb_i)+C_{11}(r_iy_j+r_jy_i),
\end{eqnarray}
where
\begin{eqnarray}
\nonumber
C_1 \!\!\!\!&:=&\!\!\!\!\ \frac{1}{2\alpha^3\Delta^2}\Big\{\Phi\alpha Q''s_0+2\alpha\Delta^2\Psi'' r_0-\Delta^2\Omega'' r_0+2\Delta^2\alpha\Omega''Qs_0\\ \!\!\!\!&&\!\!\!\!\ +4\Delta^2\alpha\Omega'Q's_0+2\alpha\Delta^2\Psi''s_0\Big\},\nonumber\\
\nonumber
C_2\!\!\!\!&:=&\!\!\!\!\ \frac{-1}{2\alpha^4\Delta^2}\Big\{2\alpha\Delta^2\Psi'' s_0-2\Omega'\Delta^2r_0+2\Omega'\Delta^2\alpha Qs_0-\Delta^2\Omega''s r_0\\ \nonumber
&& +2\Delta^2\alpha\Omega''sQs_0+4\Delta^2\alpha\Omega'Q's_0s+2\alpha\Delta^2\Psi'r_0+2\alpha\Delta^2\Psi'' s r_0\\
&& +2\alpha\Delta^2\Psi'' ss_0 +\Phi\alpha Q' s_0+\Phi\alpha Q'' s_0 s\Big\}, \nonumber\\
\nonumber
C_3\!\!\!\!&:=&\!\!\!\!\ \frac{1}{4\alpha^5\Delta^2}\Big\{4\Delta^2 s^2\Omega''\alpha Q s_0-2\Delta^2 s^2\Omega''r_0+12\alpha\Delta^2\Psi'sr_0+ 12\alpha\Delta^2\Psi's s_0\\ \nonumber
&& +4\alpha\Delta^2\Psi''s^2 r_0+4\alpha\Delta^2\Psi'' s^2 s_0+8\Delta^2 s^2 \Omega'\alpha Q's_0+2\Phi\alpha Q'' s_0s^2\\ \nonumber
&& -10\Omega'\Delta^2sr_0+12\Omega'\Delta^2 s\alpha Q s_0+6\Phi\alpha Q's_0 s-3\Phi r_0\Big\},\\
\nonumber C_4\!\!\!\!&:=&\!\!\!\! \frac{-1}{4\alpha^3\Delta^2}\Big\{4\alpha\Delta^2\Psi's s_0-\Phi r_0-2\Omega'\Delta^2s r_0+4\Omega'\Delta^2s \alpha Q s_0\\ \!\!\!\!&&\!\!\!\!\ +4\alpha\Delta^2\Psi' s r_0+2\Phi\alpha Q' s_0 s\Big\},\nonumber\\
C_5 \!\!\!\!&:=&\!\!\!\!\ \frac{-\Omega'}{\alpha^2},\ \  C_6:=\frac{2\Delta^2s\Omega'+\Phi}{2\alpha^3\Delta^2},\,\,\,C_7:=\frac{-\Phi}{2\alpha\Delta^2}, \nonumber\\
C_8 \!\!\!\!&:=&\!\!\!\!\ \frac{1}{2\alpha\Delta^2}\{2\Omega'\Delta^2 Q+2\Delta^2 \Psi'+\Phi Q'\}, \nonumber\\
C_9 \!\!\!\!&:=&\!\!\!\!\ \frac{-s}{\alpha}C_8, \ \  C_{10}:=\frac{\Psi'}{\alpha}, \,\,\,\,\, C_{11}:=\frac{-s}{\alpha}C_{10}.\nonumber
\end{eqnarray}
\end{prop}
The formula of $E$-curvature of Randers metrics and Kropina metrics computed from Proposition \ref{E} coincides with  the one computed in \cite{BY}.

\bigskip

It is easy to see that $F$ is of isotropic mean Berwald curvature if and only if
$F$ is of weak isotropic S-curvature. Hence, we consider an $(\alpha,\beta)$-metric $F=\alpha\phi(\beta/\alpha)$ with weak isotropic S-curvature, $ {\bf S}= (n+1) c F+\eta$, where $\eta=\eta_i(x)y^i$ is a 1-form on underlying manifold $M$.  Using the same method used in \cite{ChSh3}, one can obtains that the condition weak isotropic S-curvature $ {\bf S}= (n+1) c F+\eta$ is equivalent to the following equation
\be
\alpha^{-1} \frac{\Phi}{2\Delta^2} (r_{00}-2 \alpha Q s_0) - 2 \Psi (r_0+s_0)  =  -(n+1) c F +\widetilde{\theta},\label{ABC}
\ee
where
\be
\widetilde{\theta }:= -\frac{f'(b)}{bf(b)} (r_0+s_0)-\eta.
\label{theta}
\ee
To simplify the equation (\ref{ABC}), we  choose  special coordinates
$\psi:  (s,u^A) \to (y^i)$ as follows
\be y^1 = \frac{s}{\sqrt{b^2-s^2} }\bar{\alpha}, \ \ \ \ \ \ y^A = u^A,\label{specialcoordinates}
\ee
where
\[
\bar{\alpha} =\sqrt{\sum_{A=2}^n (u^A)^2}.
\]
Then
\be
\alpha =\frac{b}{\sqrt{b^2-s^2} }\bar{\alpha}, \ \ \ \ \
\beta = \frac{bs}{\sqrt{b^2-s^2} } \bar{\alpha}. \label{abara}
\ee Fix an arbitrary point $x$.
Take a local coordinate system   at  $x$ as in (\ref{specialcoordinates}). We have
\[ r_1 = b r_{11}, \ \ \ \ r_{A} = b r_{1A},\]
\[ s_1 = 0, \ \ \ \ s_{A} = b s_{1A}.\]
Let
\[ \bar{r}_{10}:= \sum_{A=2}^nr_{1A} y^{A}, \ \ \ \ \ \bar{s}_{10}:= \sum_{A=2}^ns_{1A} y^{A}\ \ \ \ \ \bar{r}_{00} := \sum_{A,B=2}^n r_{AB}y^{A} y^{B},
\]
\[
\bar{r}_0 := \sum_{A=2}^n r_{A} y^{A} \ \ \ \ \ \bar{s}_{0}:=\sum_{A=2}^n s_{A} y^{A}.
\]
Put
\[
\widetilde{\theta}= t_i y^i-\eta_iy^i.
\]
Then $t_i$ are given by
\be  t_1 = -\frac{f'(b)}{f(b)}  r_{11},
\ \ \ \ t_{A} = -\frac{f'(b)}{f(b)} (r_{1A}+s_{1A}).\label{specialti}
\ee
From (\ref{specialcoordinates}), we have
\be
r_{0}=\frac{sbr_{11}}{\sqrt{b^{2}-s^{2}}}\bar{\alpha}+b\bar{r}_{10}, \ \ \ \ s_{0}=\bar{s}_{0}=b\bar{s}_{10}, \label{r0s0}
\ee
and
\be
r_{00}=\frac{s^{2}\bar{\alpha}^{2}}{b^{2}-s^{2}}r_{11}+2\frac{s\bar{\alpha}}{\sqrt{b^{2}-s^{2}}}\bar{r}_{10}+\bar{r}_{00},\label{r0000}
\ee
\be
\widetilde{\theta}=t_1\frac{s}{\sqrt{b^2-s^2}}\bar{\alpha}-\frac{f'(b)}{f(b)}\bar{r}_{10}
-\frac{f'(b)}{f(b)}\bar{s}_{10}-\eta. \label{bar alpha theta}
\ee
Substituting (\ref{r0s0}), (\ref{r0000}) and (\ref{bar alpha theta}) into (\ref{ABC}) and by using (\ref{abara}), we find that  (\ref{ABC}) is equivalent to the following  equations:
\be
\frac{\Phi}{2\Delta^2} (b^2-s^2)\bar{r}_{00} = - \Big \{ s \Big ( \frac{s\Phi}{2\Delta^2} - 2 \Psi b^2  \Big )  r_{11}+ (n+1) c b^2 \phi - s b t_1     \Big \}\bar{\alpha}^2,  \label{C1}
\ee
\be
\Big (  \frac{s\Phi}{\Delta^2} - 2\Psi b^2  \Big ) ( r_{1A}+s_{1A})- (b^2Q+s) \frac{\Phi}{\Delta^2}  s_{1A} +b\eta_A- b t_{A}=0.\label{C2}
\ee
\be \label{C3}
\eta_1=0.
\ee
Let
\[
\Upsilon:= \Big [ \frac{s\Phi}{\Delta^2}  -2 \Psi b^2 \Big ]'.
\]
We see that
 $\Upsilon =0$ if and only if
\[ \frac{s\Phi}{\Delta^2}  - 2 \Psi b^2  = b^2 \mu,\]
where $\mu =\mu (x)$ is   independent of $s$.
\par
Let us suppose that $\Xi=\frac{(b^2 Q+s) \Phi}{\Delta^2}$ is not constant. Now we shall divide the proof into two cases: (i) $\Upsilon=0$ and (ii) $\Upsilon\not=0$.

\subsection{$\Upsilon=0$}
First, note that
 $\Upsilon=0$ implies that
\be
 \frac{s\Phi}{\Delta^2} -2 \Psi b^2 = b^2 \mu,\label{Bbb}
\ee
where $\mu =\mu (x)$ is a function on $M$ independent of $s$. First, we prove the following.

\begin{lem}\label{lem5.1}
 Let $(M, F)$ be a n-dimensional Finsler manifold. Suppose that $F=\alpha \phi (\beta/\alpha)$ be an $(\alpha,\beta)$-metric and $\Upsilon=0$. If $F$ has  weak isotropic S-curvature, ${\bf S}= (n+1)c F+\eta$, then
$\beta$ satisfies
\be
r_{ij} = k a_{ij} - \e b_i b_j + \frac{1}{b^2} (r_i b_j + r_j b_i ),\label{rrrr}
\ee
where $k=k(x)$, $\e=\e(x)$,
and
 $\phi=\phi(s)$ satisfies the following ODE:
\be
(k-\epsilon s^2) \frac{\Phi}{2\Delta^2} = \Big \{ \nu +(k -\e b^2)\mu \Big \}  s -(n+1) c \phi,\label{r4}
\ee
where $\nu=\nu (x)$.
If $s_0\not =0$, then  $\phi$ satisfies the following additional ODE:
\be
\frac{\Phi}{\Delta^2} (Q b^2+s) = b^2 (\mu + \lambda), \label{r5*}
\ee
where $\lambda=\lambda(x) $.
\end{lem}
{\it Proof}:
Since $\Phi\not=0$ and  $\bar{r}_{00}$  and $\bar{\alpha}$ are independent of $s$, it follows from (\ref{C1}) and (\ref{C2}) that
 in a special coordinate system $(s, y^a)$ at a point $x$, the following  relations hold
\be
r_{AB} = k \delta_{AB},\label{r1}
\ee
\be
s\Big ( \frac{s \Phi}{2\Delta^2} - 2 \Psi b^2    \Big ) r_{11} +(n+1) c b^2 \phi + k \frac{\Phi}{2\Delta^2} (b^2-s^2) = b s t_1, \label{r2}
\ee
\be
\Big ( \frac{s\Phi}{\Delta^2} -2 \Psi b^2 \Big )( r_{1A} +s_{1A}) - (b^2Q+s) \frac{\Phi}{\Delta^2}  s_{1A} - b t_{A} =-b\eta_A, \label{r3}
\ee
where $k=k(x)$ is  independent of $s$. Let
\[ r_{11}= -(k-\epsilon b^2).\]
Then (\ref{rrrr}) holds.
By (\ref{Bbb}), we have
\[ \frac{s\Phi}{2\Delta^2} - 2 \Psi b^2
= b^2 \mu -\frac{s\Phi}{2\Delta^2}.\]
Then
(\ref{r2}) and (\ref{r3}) become
\be
b (k-\e s^2)\frac{\Phi}{2\Delta^2} = s t_1 + sb\mu (k-b^2\e) -(n+1) c b \phi.\label{r2*}
\ee
\be
b^2 \mu (r_{1A}+s_{1A})  -\frac{\Phi}{\Delta^2} (Q b^2+s)  s_{1A} - b t_{A} =-b\eta_A.\label{r3*}
\ee
Letting  $t_1 = b \nu$ in (\ref{r2*}) we get (\ref{r4}). Now, suppose that  $s_0\not=0$.
Rewrite (\ref{r3*}) as
\[  \Big \{  b^2 \mu - \frac{\Phi}{\Delta^2} (Q b^2+s) \Big \} s_{1A}
= b t_A -b\eta_A- b^2 \mu r_{1A}.\]
We can see that there is a function $\lambda =\lambda (x)$  on $M$ such that
\[
 \mu b^2 - \frac{\Phi}{\Delta^2} (Q b^2+s) = -b^2 \lambda.
\]
This gives (\ref{r5*}).\qed

\bigskip

\begin{lem} {\rm(\cite{ChSh3})} \label{lemb=c2}  Let $F =\alpha \phi (\beta/\alpha)$ be an $(\alpha,\beta)$-metric. Assume that $$\phi \not= k_1 \sqrt{1+k_2  s^2}+k_{3}s$$ for any constants $k_1>0, k_{2}$ and $k_3$.
If $\Upsilon=0$, then $b = constant$.
\end{lem}

An $(\alpha, \beta)$-metric is called Randers-type if $\phi = k_1 \sqrt{1+k_2  s^2}+k_{3}s$ for any constants $k_1>0, k_{2}$ and $k_3$. Now, we consider the equivalency of the notions  weak isotropic S-curvature and   isotropic S-curvature for a non-Randers type $(\alpha,\beta)$-metric.

\begin{lem}\label{propPart2}
Let $F =\alpha \phi (\beta/\alpha)$ be a non-Randers type $(\alpha,\beta)$-metric.
Suppose that $\Xi$ is not constant and $\Upsilon=0 $.
Then $F$ is of weak isotropic S-curvature if and only if $F$ is of isotropic S-curvature.
\end{lem}
{\it Proof}:  It is sufficient to prove that if  $F$ is of weak isotropic S-curvature, then $F$ is of isotropic S-curvature. By $db=(r_0+s_0)/b$ and  Lemma \ref{lemb=c2}, we have
\[ r_0 + s_0 =0.\]
Then by the formula of S-curvature of an $(\alpha,\beta)$-curvature, we get
\[ {\bf S} = -\alpha^{-1} \frac{\Phi}{2\Delta^2} \Big \{ r_{00} - 2 \alpha Q s_0 \Big \} .\]
By Lemma \ref{lem5.1},
\[ r_{00} = (k -\e s^2)\alpha^2 +\frac{2s}{b^2} r_0 \alpha.\]
Then
\[
 {\bf S}=- (k-\e s^2) \frac{\Phi}{2\Delta^2}  \alpha
+  \frac{\Phi}{b^2\Delta^2} (b^2Q+s) s_0 .
\]
By (\ref{r4}), we have
\be
{\bf S} = - s \Big \{ \nu + (k-\e b^2) \mu \Big \} \alpha
+ \frac{\Phi}{b^2\Delta^2} (b^2Q+s) s_0 + (n+1) c \phi \alpha. \label{mmm}
\ee
Since ${\bf S}= (n+1) c F+\eta$, then by (\ref{mmm}) we obtain the following
\be
- s \Big \{ \nu + (k-\e b^2) \mu \Big \} \alpha
+ \frac{\Phi}{b^2\Delta^2} (b^2Q+s) s_0=\eta. \label{mmm*}
\ee
Letting $ y^i = \delta b^i$ for a sufficiently small $\delta >0$ yields
\[  - \delta \Big \{ \nu + (k-\e b^2) \mu \Big \} b^2 =\delta \eta_ib^i.\]
It is easy to see that in the special coordinate $\eta_ib^i=0$, hence in general $\eta_ib^i=0$.
We conclude that
\be
 \nu + (k-\e b^2) \mu =0.\label{d1}
\ee
Then (\ref{mmm*}) reduce to
\be
\frac{\Xi}{b^2} s_0=\eta.\label{aki}
\ee
If $s_0\not=0$, then from the last equation, we obtain that $\Xi$ is constant, which is excluded here.
Hence, we have $s_0 =0$. Thus by (\ref{aki}), we conclude that $\eta=0$ and $F$ has isotropic S-curvature ${\bf S}=(n+1)cF$.
\qed

\subsection{$\Upsilon\not=0$}
Here, we  consider the case when $\phi=\phi(s)$ satisfies
\be
\Upsilon \not=0\label{CB}
\ee
First we need the following
\begin{lem}\label{thm4.1}
Let $F=\alpha\phi(s), s=\beta/\alpha$,  be an $(\alpha,\beta)$-metric on an $n$-dimensional manifold. Assume that $\Upsilon \not=0$.
Suppose that $F$ has  weak isotropic S-curvature, ${\bf S}= (n+1)c F+\eta$.
 Then
\be
r_{ij} = k a_{ij}-\e b_i b_j -\lambda (s_i b_j + s_j b_i)  ,\label{SS1}
\ee
where $\lambda=\lambda(x), k=k(x) $ and $\e=\e(x)$ are scalar functions of $x$ and
\be
  -2 s(k-\e b^2) \Psi + (k-\e s^2) \frac{\Phi}{2\Delta^2}  +(n+1) c \phi -s \nu=0,\label{SS2}
\ee
where
\be
\nu:= - \frac{f'(b)}{bf(b)} (k-\e b^2). \label{SS2a}
\ee
If in addition $s_0\not=0$, i.e., $s_{A_o} \not =0$ for some $A_o$, then
\be -2 \Psi -\frac{Q\Phi}{\Delta^2} -\lambda \Big ( \frac{s\Phi}{\Delta^2} -2\Psi b^2  \Big )=\delta, \label{SS3}
\ee
where
\be
 \delta:= - \frac{ f'(b)}{b f(b)} (1-\lambda b^2)-\frac{\eta_{A_o}}{s_{A_o}}.\label{S4}
\ee
\end{lem}
{\it Proof}:
By assumption,  $ \Phi \not\equiv 0$. Similar to the proof of Lemma \ref{lem5.1},
it follows from (\ref{C1}) that
there is a function $k=k(x)$  independent of $s$,   such that
\be
\bar{r}_{00} = k \bar{\alpha}^2,\label{C1a}
\ee
\be
 s\Big ( \frac{s\Phi}{2\Delta^2} - 2 \Psi b^2 \Big) r_{11}+ (n+1) c b^2 \phi   +  k \frac{\Phi}{2\Delta^2} (b^2-s^2) =sbt_1.  \label{C1b}
\ee
Let
\[ r_{11}= k -\e b^2,\]
where $\e =\e (x)$ is  independent of $s$.
By (\ref{specialti}),
\[
t_1 = b\nu,
\]
where $\nu$ is given by (\ref{SS2a}). Plugging them into (\ref{C1b}) yields (\ref{SS2}).

\bigskip
Suppose that $s_0=0$. Then
\[ b s_{1A}= s_{A}=0.\]
Then (\ref{C2}) is reduced to
\be
\Big ( \frac{s\Phi}{\Delta^2}-2\Psi b^2  \Big ) r_{1 A} - b t_{A}=-b\eta_{A}.\label{df}
\ee
By assumption, $\Upsilon \not=0$, we know that $ \frac{s\Phi}{\Delta^2} -2 \Psi b^2$ is not independent of $s$.  It follows from (\ref{df}) that
\[ r_{1A}=0, \ \ \ \ \ \ t_{A}=\eta_A.\]
The above identities together with $r_{11} = k-\e b^2$ and $t_1 = b \nu$ imply the following identities
\be
r_{ij}  = k a_{ij} -\e b_i b_j.
\ee
Now, suppose that $s_0\not=0$. Then
\[
s_{A_o} = b s_{1A_o}\not=0
\]
for some $A_o$. Differentiating (\ref{C2}) with respect to $s$ yields
\be
\Upsilon\;  r_{1A} - \Big [ \frac{Q\Phi}{\Delta^2} + 2 \Psi     \Big ]' b^2 s_{1A} =0.\label{C2**}
\ee
Let
\[ \lambda := -\frac{ r_{1A_o}}{b^2 s_{1A_o} }.\]
Plugging it into (\ref{C2**}) yields
\be
 -\lambda \Upsilon - \Big [ \frac{Q\Phi}{\Delta^2} + 2 \Psi     \Big ]' =0. \label{C2***}
\ee
It follows from (\ref{C2***}) that
\[ \delta:= - \frac{Q\Phi}{\Delta^2}-2\Psi -\lambda \Big [ \frac{s\Phi}{\Delta^2}- 2\Psi b^2   \Big ] \]
is a number independent of $s$.
By assumption that $\Upsilon\not=0$, we obtain
from (\ref{C2**}) and (\ref{C2***}) that
\be
r_{1A}+\lambda b^2 s_{1A}=0.\label{CC2}
\ee
(\ref{C1a}) and (\ref{CC2}) together with $ r_{11} =k -\e b^2$ imply that
\be
r_{ij} + \lambda ( b_i s_j + b_j s_i) =  k a_{ij} -\e b_i b_j.
\ee
By (\ref{specialti}) and (\ref{CC2}),
\[ t_{A} = \frac{f'(b)}{f(b)} (b^2\lambda-1) s_{1A}.\]
On the other hand,
by (\ref{C2}) and (\ref{CC2}), we obtain
\[ b t_{A} = \delta b^2 s_{1A}+b\eta_A.\]
Combining the above identities, we get (\ref{S4}).
\qed

\bigskip

\begin{lem}\label{lem5.2} Let $F=\alpha \phi(s)$, $s=\beta/\alpha$, be an $(\alpha,\beta)$-metric.
Suppose that  $\phi=\phi(s)$ satisfies (\ref{CB}) and $\phi\not= k_1 \sqrt{1+k_2 s^2} +k_3 s$ for any constants $k_1>0$,  $k_2$  and $k_3$. If $F$ has weak isotropic S-curvature, then
\[ r_j+s_j=0.\]
\end{lem}
{\it Proof}:  Suppose that $r_j +s_j\not=0$, then $db=(r_0+s_0)/b\neq 0$ and hence
$b:=\|\beta_x\|_{\alpha} \not= constant$ in a neighborhood. We view $b$ as a variable in (\ref{SS2}) and (\ref{SS3}).  Since $\phi =\phi(s)$ is a function independent of $x$, (\ref{SS2}) and (\ref{SS3}) actually give rise infinitely many ODEs on $\phi$. First, we consider (\ref{SS2}).
Let
\[ eq:= \Delta^2 \Big \{ -2s(k-\e b^2)\Psi + (k-\e s^2) \frac{\Phi}{2\Delta^2}+(n+1) c \phi -s \nu\Big \}.\]
We have
\[ eq= \Xi_0+\Xi_2 b^2 + \Xi_4 b^4 ,\]
where $\Xi_0, \Xi_2, \Xi_4$ are independent of $b$ and
\[ \Xi_4 := \Big \{  (\e -\nu) s + (n+1) c \phi \Big \} \frac{\phi^2}{(\phi-s \phi')^4} (\phi'')^2.
\]
It follows from (\ref{SS2}) that $eq=0$. Thus
\[ \Xi_0=0, \ \ \ \ \Xi_2=0, \ \ \ \ \Xi_4=0.\]
Since $\phi''\not=0$, the equation
$\Xi_4=0$ is equivalent to the following ODE
\[ (\e -\nu) s + (n+1) c \phi  =0.\]
We conclude that
\[ \e =\nu, \ \ \ \ \ \ c=0.\]
Then by a direct computation we get
\[  \Xi_0+\Xi_2 s^2 =
-\frac{1}{2} (1+sQ) \Big \{ (n-1)(k-\e s^2) (Q-sQ') + 2k Q + 2\e s      \Big \}.\]
By  $\Xi_0=0$ and $\Xi_2=0$ it follows that
\be
(n-1)(k-\e s^2) (Q-sQ') + 2k Q + 2\e s   =0,\label{eqQQQ}
\ee
Suppose that $(k, \e)\not=0$. We claim that $k\not=0$. If this is not true, i.e.,  $k=0$, then $\e\not=0$ and (\ref{eqQQQ}) is reduced to
\[
-(n-1) s (Q-sQ')+2  =0.
\]
Letting $s=0$, we get a contradiction.

Now we have  $k\not=0$.  It is easy to see that $Q(0)=0$. Let
\[ \tilde{Q}:=Q(s)-sQ'(0).\]
Plugging it into  (\ref{eqQQQ}) yields
\be
(n-1)(k-\e s^2) (\tilde{Q}-s\tilde{Q}') + 2k \tilde{Q}
+ 2 ( k Q'(0)+\e) s=0.\label{Q}
\ee
Differentiating  of the equation (\ref{Q}) with respect to $s$ and then letting $s=0$,   yields  
\[
kQ'(0)+\e =0.
\]
Then (\ref{Q}) is reduced to following
\[ (n-1)(k-\e s^2) (\tilde{Q}-s\tilde{Q}') + 2k \tilde{Q}=0.\]
Solving above ODE, we obtain
\[ \tilde{Q} = c_1 s\,e^{\frac{2k}{n-1}\int_0^s{\frac{1}{u(k-\epsilon u^2)}du}}.\]
Since $\tilde{Q}'(0)=0$, we have $c_1=0$. Hence, $\tilde{Q}=0$. We get
\[  Q(s)-sQ'(0)=0.\]
Then it follows that
\[ Q(s) = Q'(0)s.\]
In this case, $\phi= c_1 \sqrt{1+c_2s^2}$ where $c_1>0$ and $c_2$ are numbers independent of $s$. This case is excluded in the assumption.
Therefore $k=0$ and $\e=0$. Then  (\ref{SS1}) is reduced to
\[
r_{ij} =-\lambda (s_j b_i + s_i b_j).
\]
Then
\[ r_j + s_j = (1-\lambda b^2) s_j.\]
At the beginning of the proof, we suppose that $r_j+s_j \not=0$. Thus, we  conclude that
\[
1-\lambda b^2\not=0,\ \  \textrm{and}\ \   s_j \not=0.
\]
By Lemma \ref{thm4.1},  $\phi=\phi(s)$ satisfies (\ref{SS3}).
Let
\[
 EQ:=\Delta^2 \Big \{ -2 \Psi - \frac{Q\Phi}{\Delta^2} -\lambda \Big ( \frac{s\Phi}{\Delta^2} - 2 \Psi b^2 \Big )-\delta\Big \}.
\]
We have
\[ EQ = \Omega_0 + \Omega_2 b^2 + \Omega_4 b^4,\]
where $\Omega_0, \Omega_2, \Omega_4$ are independent of $b$ and
\[ \Omega_4 =  (Q')^2 (\lambda -\delta).\]
By (\ref{SS3}), $EQ=0$. Thus
\[\Omega_0=0,\ \ \ \Omega_2=0, \ \ \  \Omega_4=0.\]
Since $Q'\not=0$, $\Omega_4=0$ implies that
\[\delta =\lambda.\]
By a direct computation, we get
\[ \Omega_0+\Omega_2 s^2 = (1+sQ) \Big \{  (n+1) Q (Q-sQ')-Q' +\lambda \Big [n s (Q-sQ')-1 \Big ] \Big \}.\]
The equations $\Omega_0=0$ and $\Omega_2=0$ imply that 
\[
\Omega_0+\Omega_2 s^2=0,
\]
that is
\[  (n+1) Q (Q-sQ')-Q' +\lambda \Big [n s (Q-sQ')-1 \Big ]  =0.\]
We obtain
\[ Q=
-\frac{[ k_{0}n(n+1)-1] \lambda s\pm \sqrt{\lambda k_{0}  (k_{0} (1+n)^2 -1+\lambda s^2)}}{k_{0}(n+1)^2-1}.
\]
Plugging it into $\Omega_2=0$ yields
\[ k_{0} \lambda =0.\]
Then
\[ Q= \frac{\lambda s}{ k_{0}(n+1)^2-1}.\]
This implies that 
\[
\phi = k_1 \sqrt{1+k_2 s^2},
\]
where $k_1>0$ and $k_2$ are numbers independent of $s$.
This case is excluded in the assumption of the lemma. Therefore, $r_j+s_j =0$.
\qed

\bigskip
\begin{prop}\label{propA}
Let $F=\alpha \phi(s)$, $s=\beta/\alpha$, be an $(\alpha,\beta)$-metric.  Suppose that  $\phi =\phi(s)$ satisfies  (\ref{CB}) and $\phi\not= k_1 \sqrt{1+k_2 s^2} +k_3 s$ for any constants $k_1>0$,  $k_2$  and $k_3$. Suppose that $\Xi$ is not constant.  If $F$ is of weak isotropic S-curvature,
${\bf S}= (n+1)c F+\eta$,
 then
\be
r_{ij}= \e (b^2 a_{ij}-b_ib_j), \ \ \ \ \ s_j =0,\label{S1AAA}
\ee
where  $\e=\e(x)$ is a scalar function on $M$
and $\phi =\phi(s)$ satisfies
\be
\e (b^2- s^2) \frac{\Phi}{2\Delta^2} = -(n+1) c \phi.\label{S2AAA}
\ee
\end{prop}
{\it Proof}:
Contracting (\ref{SS1}) with $b^i$ yields
\be
r_j + s_j = (k-\e b^2) b_j + (1-\lambda b^2) s_j. \label{SS1***}
\ee
By Lemma \ref{lem5.2}, $r_j + s_j =0$.
It follows from (\ref{SS1***}) that
\be
 (1-\lambda b^2) s_j + (k-\epsilon b^2) b_j =0.\label{srjj}
\ee
Contracting (\ref{srjj}) with $b^j$ yields
\[
(k-\epsilon b^2 ) b^2 =0.
\]
We get
\[ k = \epsilon b^2.\]
Then (\ref{SS1}) is reduced to
\[ r_{ij} = \e ( b^2 a_{ij}- b_i b_j) -\lambda (b_i s_j +b_j s_i ).\]
By (\ref{SS2a}),
\[ \nu=0.\]
Then (\ref{SS2}) is reduced to (\ref{S2AAA}).

We claim that $s_0=0$.
Suppose that $s_0\not=0$.
By (\ref{srjj}), we conclude that
\[
 \lambda = \frac{1}{b^2}.
\]
By (\ref{S4}),
\[
\delta =-\frac{\eta_{A_o}}{s_{A_o}}.
\]
It follows from (\ref{SS3}) that
\[
\frac{(b^2 Q+s) \Phi }{\Delta^2}= \frac{b\eta_{A_o}}{s_{A_o}},
\]
which implies that $\Xi$ is constant. This is impossible by the assumption on non-constancy of $\Xi$. Therefore, $s_j=0$. This completes
the proof.
\qed

\bigskip

By Proposition \ref{propA} and Lemma \ref{Cheng-Shen}, we have the following.
\begin{cor}
Let $F=\alpha \phi(s)$, $s=\beta/\alpha$, be a non-Randers type  $(\alpha,\beta)$-metric.  Suppose that $\Upsilon\neq 0$ and $\Xi$ is not constant.  Then $F$ is of weak isotropic S-curvature, if and only if it is of isotropic S-curvature.
\end{cor}

\bigskip

\noindent
Behzad Najafi\\
Department of Mathematics and Computer Sciences\\
Amirkabir University\\
Tehran. Iran\\
Email:\ behzad.najafi@aut.ac.ir
\bigskip

\noindent
Akbar Tayebi\\
Faculty  of Science, Department of Mathematics\\
University of Qom \\
Qom. Iran\\
Email: akbar.tayebi@gmail.com
\bigskip

\end{document}